\documentclass[12pt]{article}
\usepackage{amsmath}

\usepackage{amsthm} %\usepackage{thmdefs}
\theoremstyle{plain} %% This is the default
\newtheorem{theorem}{Theorem}[section]

\newtheorem{lemma}[theorem]{Lemma}

\numberwithin{equation}{section}

\newcommand{\al}{\alpha}

\newcommand{\br}{\mathbf {R}}
\newcommand{\om}{\omega}
\newcommand{\ot}{\otimes}
\newcommand{\ld}{\ldots}
\newcommand{\p} {\partial}
\newcommand{\g} {\mathsf g }
\newcommand{\nab}{\nabla}

\begin{document}
\title{Leibniz Cohomology and the Calculus of Variations}
\author{Jerry M. Lodder}
\date{Math Sciences, Dept. MB  \\ New Mexico State University \\
      Las Cruces, NM  88003}
\maketitle

\section{Introduction}

In a previous publication \cite{Lodder} a definition of Loday's
Leibniz cohomology, $HL^*$, \cite{LP} was proposed for differentiable
manifolds.  In particular every $k$-tensor $\om$ (from classical
differential geometry) is a cochain in the Leibniz complex.  Although
the Leibniz coboundary, $d \om$, is not necessarily a $(k+1)$-tensor,
$d \om$ remains a local operator on vector fields with the value of $d
\om$ at a point $p$ in the manifold $M$ determined by the 
values of $\om$ in an arbitrary open neighborhood containing $p$.
With this writing we offer an explicit formula for $d \om$ in a local
coordinate chart, and provide a geometric interpretation of $d \om$ in
terms of the calculus of variations. If $\om$ is the metric two-tensor
on a Riemannian manifold, then the local expression for $d \om$
involves the Christoffel symbols, while the global definition of $d
\om$ reduces to the first variation formula for arc length.  More
generally the Leibniz coboundary of any two tensor $\om$ can be
written in terms of the necessary conditions to achieve a minimum (or
maximum) value of $\int \om$ over a locally immersed curve or
surface.  The paper closes with the computation of the Leibniz
coboundary of the Riemann curvature tensor $R$ in terms of its
covariant derivative $\nabla R$.

Section two of the paper begins with a brief recollection of $HL^*$
for a differentiable manifold $M$, and proceeds with the foundational
material needed to prove that $d \om$ is a local operator.  Section
three contains the results for two tensors and the calculus of
variations.  The final section provides the local coboundary formula
for arbitrary $k$-tensors as well as the global coboundary of the
Riemann curvature tensor.  For more background material about
$HL^*(M)$ and in particular calculations of $HL^*$ for Euclidean
$n$-space (which are highly non-trivial), 
see \cite{Lodder}.  For more information about Leibniz
homology and cohomology, see \cite{Ens.Math.}
\cite{Bourbaki} \cite{LP}.

\section{The Leibniz Coboundary as a Local \\ Operator}

We begin by reviewing the definition of Leibniz cohomology for
differentiable manifolds \cite{Lodder}, and show that the Leibniz
coboundary of a $k$-tensor, $d \om$, is a local operator, i.e. $d \om$
at $p \in M$ is determined by the value of $\om$ on an arbitrary open
neighborhood of $p \in M$.  This permits (in later sections) the
formulation of $d \om$ in terms of a local coordinate chart.  Let $M$
be a differentiable ($C^{\infty}$) manifold of dimension $n$, $\chi
(M)$ the Lie algebra of $C^{\infty}$ vector fields on $M$, and
$C^{\infty}(M)$ the algebra of $C^{\infty}$ real-valued functions $f:M
\to \br$.  Recall that $C^{\infty}(M)$ is a left representation of
$\chi (M)$ via 
\begin{align*}
 & \chi (M) \, {\underset {\br} {\otimes}} \, C^{\infty}(M) \to  C^{\infty}(M)\\
 & [X, \, f] \mapsto X(f),   
\end{align*}
where $X(f)$ is the Lie derivative of $f \in C^{\infty}(M)$ in the
direction $ X \in \chi (M)$.  Let
$$ C^k (M) = \text{Hom}_{\br}^{c} ( \chi(M)^{\otimes k}, \,
	C^{\infty}(M)),  \ \ \ k \geq 0 , $$
denote the $\br$-vector space of continuous homomorphisms
$$ \al :  \chi(M)^{\otimes k} \to C^{\infty}(M) $$
in the strong $C^{\infty}$ topology.  See \cite{Hirsch} for a discussion of
this topology.  Then the Leibniz cohomology of $M$ with coefficients
in $C^{\infty}(M)$, written 
$$  HL^* ( \chi(M); \, C^{\infty}(M)),  $$
is the homology of the cochain complex
$$  C^0 (M) \to C^1 (M) \to \ \ldots \ \to C^k (M) \ {\overset {d}{\to}}
	\ C^{k+1} (M) \to \ \ldots \, , $$
where
\begin{equation} \label{2.1}
\begin{split}
  & d \al  (X_1 \ot  X_2 \ot \, \ldots \, \ot X_{k+1}) = \\
  & \sum_{i=1}^{k+1} (-1)^{i+1} X_i \big( \al (X_1 \ot \, \ldots \,
   \hat{X}_i \,  \ldots \ot X_{k+1}) \big) + \\
  & \sum_{1 \leq i < j \leq k+1} (-1)^{j+1} \al \big( X_1 \ot \, \ldots
   \, \ot X_{i-1}\ot   [X_i, \, X_j]\ot X_{i+1}\ot \, \\  
  & \hskip2.5in \ldots \, 
	\hat{X}_j \, \ldots \,\ot X_{k+1}\big) .  
\end{split}
\end{equation}

Let $\om$ be a $k$-tensor on $M$, i.e.
$$  \om : M \to T^* (M)^{\ot k}  $$
is a $C^{\infty}$ section of the $k$-fold tensor product of the
cotangent bundle.  Then $\om$ determines an element of 
$$  {\text{Hom}}^{c}_{\br} ( \chi(M)^{\ot k}, \, C^{\infty}(M))  $$
via $\om (X_1 \ot X_2 \ot \, \ldots \, \ot X_k ) : M \to \br$
$$ \om (X_1 \ot X_2 \ot \, \ldots \, \ot X_k )(p) = \om ( X_1 (p) \ot
	X_2 (p) \ot \, \ldots \, X_k (p) ).  $$
Although $d \om$ is not necessarily a $(k+1)$-tensor \cite{Lodder},
the following local result remains valid.
\begin{lemma} \label{2.2}
Let $\om$ be a $k$-tensor on $M$ and $O \subset M$ open. If 
$$ X_1, \ X_2, \ \ldots, \ X_{k+1}, \ Y_1, \ \ldots, \ Y_{k+1} \in \chi (M)$$
with $X_i = Y_i$, $i = 1, \, 2, \, \ldots \, k+1$, on $O$, then
$$ d \om (X_1 \ot \, \ld \, \ot X_{k+1})(p) = 
   d \om (Y_1 \ot \, \ld \, \ot Y_{k+1})(p)   $$
for all $p \in O$.
\end{lemma}
\begin{proof}
Since the Lie bracket at $p$
$$  [X_i , \, X_j ](p)  $$
is determined by the values of $X_i$ and $X_j$ on an open set
containing $p$, we have
$$  [X_i , X_j] (p) = [Y_i , \, Y_j](p), \ \ \ p \in O . $$
Also, the (Lie) derivative of a function $f: M \to \br$ at $p$ is
determined by the values of $f$ on an open set containing $p$.  Thus, 
\begin{align*}
 X_i \big( \om ( & X_1 \ot \, \ldots \,
   \hat{X}_i \,  \ldots \ot X_{k+1}) \big) (p) =  \\
 & Y_i \big( \om ( Y_1 \ot \, \ldots \,
   \hat{Y}_i \,  \ldots \ot Y_{k+1}) \big) (p) , \ \ \ p \in O .
\end{align*}
\end{proof}

Let $x : U \to \br^n$ be a coordinate chart for $M$, $p \in U$ fixed
with $x(p) = 0 \in \br^n$.  If $X_1, \ \ld , \ X_{k+1}$ are
$C^{\infty}$ vector fields on $U$, we may define
$$  d\om (X_1 \, \ot \ \ld \ \ot \, X_{k+1})(p)  $$
for a $k$-tensor as follows.  Let $g : \br^n \to \br$ be a
$C^{\infty}$ function with
$$  g(v) = \begin{cases} 1, & ||v|| \leq 1 , \\
			 0, & ||v|| \geq 2 .\end{cases}  $$
Then $g \circ x : U \to \br$ is $C^{\infty}$ and may be extended to a
$C^{\infty}$ function $\varphi : M \to \br$ via
$$  \varphi(q) = \begin{cases} 0, & q \in M - U , \\
			      (g \circ x)(q), & q \in U .\end{cases} $$
Define $C^{\infty}$ vector fields $Y_i$ on $M$ by setting
$Y_i = \varphi \, X_i$ on $U$ and $Y_i = 0$ on $M - U$.  Clearly
$$  Y_i = X_i \ \ \ {\text{on}}\ \ \ 
	O = x^{-1}\Big( \, \{ \, v \in \br^n \ | \  
	\parallel v \parallel < 1 \, \} \, \Big).  $$
Set
\begin{equation} \label{2.3}
 d \om ( X_1  \ot \ \ld \ \ot  X_{k+1}) (p) :=
	d \om ( Y_1  \ot \ \ld \ \ot Y_{k+1}) (p) .
\end{equation}
By lemma \eqref{2.2}, the value of $d \om ( X_1 \, \ot \ \ld \ \ot \,
X_{k+1})(p)$ is independent of the choice of $g : \br^n \to \br$.  The
formula in equation \eqref{2.3} is useful for the construction of the Leibinz
coboundary of a tensor in local coordinate chart.

\section{Two Tensors and the Calculus of Variations}

In this section we compute the Leibniz coboundary of a two tensor in
terms of the local coordinate chart $(x, \ U)$, where $U \subset M$ is
open, and 
$$  x : U \to \br^n  $$ 
is a homeomorphism belonging to the atlas of charts for the
differentiable structure of $M$.  The coefficients of this coboundary
can be identified with those which occur in the optimization process
for the integral of a two tensor over an immersed curve or surface
(within $U$).  For example, the Leibniz coboundary of the metric
two-tensor (for $M$ Riemannian) can be expressed in terms of the
Christoffel symbols.  

For completeness we begin with a one-form (i.e. a one-tensor), which
has a local expression on $U$ as
$$  \om = \sum_{i=1}^{n} a_i \, dx^i ,  $$
where $a_i : U \to \br$ are $C^{\infty}$ functions.  From equation
\eqref{2.1}, the Leibniz coboundary of $\om$ agrees with the de Rham
coboundary of $\om$, and in local coordinates
\begin{equation} \label{3.1}
\begin{split}
 d \om & = \sum_{i, \, j = 1}^{n} \frac{\partial a_i}{\partial x^j} 
      \, dx^j \wedge dx^i  \\
  & = \sum_{j < i}\Big( \frac{\partial a_i}{\partial x^j}
	- \frac{\partial a_j}{\partial x^i} \Big) \, dx^j \wedge dx^i .
\end{split}
\end{equation}

We now discuss in what sense the functions 
$$  \frac{\partial a_i}{\partial x^j} - \frac{\partial a_j}{\p x^i}$$
arise from the calculus of variations.  Let 
$$  \gamma : [0, \, 1] \to U, \ \ \ \gamma(0) = p, \ \gamma(1) = q  $$
be a $C^{\infty}$ curve with a $C^{\infty}$ variation
$$  \al : ( - \epsilon , \, \epsilon) \times [0, \, 1] \to U  $$
satisfying
\begin{align*}
&  \al (0, \, t) = \gamma (t)  \\ 
&  \al (s, \, 0 ) = p \ \ \ {\text{for}} \ \  -\epsilon < s < \epsilon \\
&  \al (s, \, 1) = q \ \ \ {\text{for}} \ \  -\epsilon < s < \epsilon .
\end{align*}
We wish to investigate to what extent
$$  J(\gamma) = \int_0^1 \om (\textstyle{\frac{d \gamma}{d t}}) \, dt$$
is an extreme value (as a function of $s$) of
$$  J( \al (s)) = \int_0^1
  \om \Big( \frac{ \p \al}{\p t} (s, \, t) \Big) \, dt .  $$
Let $\gamma^i = x^i (\gamma (t)) \in \br$ be the i-th component of the
curve $\gamma (t)$.  Then
\begin{align*}
& \frac{d \gamma}{d t} = \sum_{i = 1}^n \frac{d
 \gamma^i}{d t}\, \frac{\p}{\p x^i} \, \Bigr\rvert_{\gamma(t)} \\
& \om \Big( \frac{d \gamma}{d t} \Big) = \sum_{i =1}^n a_i (\gamma
  (t)) \, \frac{ d \gamma^i}{d t} \\
& \int_0^1 \om \Big( \frac{d \gamma}{d t} \Big) \, dt =
  \int_0^1 \Big( \sum_{i =1}^n a_i (\gamma
  (t)) \, \frac{ d \gamma^i}{d t} \Big) \, dt  
\end{align*}

Recalling the treatment of the calculus of variations, \cite[p.\
438]{SpivakI}, we define
$$  F: \br^n \times \br^n \to \br  $$
by $F(x, \, y) = \sum_{i=1}^n a_i(x) \, y^i$, where $y = (y^1, \,
y^2, \, \ld \, y^n)$.  In our case
$$  a_i(x) := (a_i \, x^{-1})( x( \gamma (t))) \ \ \ \text{and} \ \ \
  y^i = \frac{d \gamma^i}{d t} . $$
A necessary condition that $J(\gamma)$ be an extreme value is that the
``Euler-Lagrange'' equations hold \cite[p.\ 438]{SpivakI}
\begin{equation} \label{3.2}
 \frac{ \p F}{\p x^{\ell}} = \frac{d}{dt} \Big( \frac{\p F}{\p y^{\ell}}\Big),
	\ \ \ell = 1, \, 2, \, \ld \, n .  
\end{equation}
Now,
\begin{align*}
& \frac{\p F}{\p x^{\ell}} \Big( \gamma(t), \, \frac{d \gamma}{d t}\Big) =
  \sum_{i=1}^n \frac{\p a_i}{\p x^{\ell}}\Big( \gamma (t) \Big) \, 
	\frac{ d \gamma^i}{d t}  \\
& \frac{\p F}{\p y^{\ell}} \Big( \gamma(t), \, \frac{d \gamma}{d t}\Big) =
  a_{\ell} (\gamma (t))  \\
& \frac{d}{dt} \Big( a_{\ell} ( \gamma (t)) \Big) = 
  \sum_{j=1}^n \frac{\p a_{\ell}}{\p x^j}\Big( \gamma (t) \Big) \, 
	\frac{ d \gamma^j}{d t}  
\end{align*}
The condition in equation \eqref{3.2} may be rewritten as 
\begin{equation} \label{3.3}
\sum_{i=1}^n \Big( \frac{\p a_i}{\p x^{\ell}} - 
  \frac{\p a_{\ell}}{\p x^i}\Big)
  (\gamma (t)) \; \frac{d \gamma^i}{d t} = 0
\end{equation}
for each $\ell = 1, \ 2, \ 3, \ \ld , \ n$.  
The coefficients arising in \eqref{3.3} also appear in \eqref{3.1}.  
Of course,
\begin{align*}
 d \om & \Bigg( \sum_{i = 1}^n \frac{d \gamma^i}{dt} \Big(
  \frac{\p}{\p x^{\ell}}\Bigr\rvert_{\gamma (t)} 
	\ot \frac{\p}{\p x^i}\Bigr\rvert_{\gamma (t)}
  \Big) \Bigg)  \\
 & = \sum_{i=1}^n \big( \frac{\p a_i}{\p x^{\ell}} - 
	\frac{\p a_{\ell}}{\p x^i}\Big)
   (\gamma(t)) \, \frac{d \gamma^i}{dt}.
\end{align*}
We now prove an identical result for the Leibniz coboundary of a
two-tensor (which is not necessarily a two-form).

Let $\om$ be a two-tensor on $M$ with local expression
$$  \om = \sum_{i, \, j = 1}^n a_{ij} \, dx^i \ot dx^j  $$
on $U$.  Let $\gamma : I^2 \to U$ be an immersion ($C^{\infty}$ is
sufficient), where
$$  I^2 = \{ \, (t_1, \, t_2) \in \br^2 \ | \ 0 \leq t_1 \leq 1, \ 0 \leq t_2
 	\leq 1 \, \}.  $$
Although the value of the integral
\begin{equation} \label{3.4}  
  J(\gamma) = \int_0^1 \int_0^1 \om \Big( \frac{\p \gamma}{\p t_1} \ot
  \frac{ \p \gamma}{\p t_2}\Big) \, dt_1 \, dt_2
\end{equation}
generally depends on the parameterization $\gamma$ (and not just the
image of $\gamma$), necessary conditions for an extreme value of
$J(\gamma)$ can still be sought.  Consider the $C^{\infty}$ variation
$$ \al : ( - \epsilon , \, \epsilon) \times I^2 \to U  $$
satisfying
\begin{align*}
& \al(0, \, t_1, \, t_2) = \gamma (t_1, \, t_2) \\
& \al (s, \, 1, \, t_2) = \gamma(1, \, t_2), \ \ \
  \al (s, \, 0, \, t_2) = \gamma(0, \, t_2), \ \ \ - \epsilon < s < \epsilon \\
& \al (s, \, t_1, \, 1) = \gamma(t_1, \, 1), \ \ \
  \al (s, \, t_1, \, 0) = \gamma(t_1, \, 0), \ \ \ - \epsilon < s <
	\epsilon . \\
\end{align*}
Then as a function of $s$, 
\begin{equation} \label{3.5}
\begin{split}
& J(\al (s)) =  \\ 
& \int_0^1 \int_0^1 \om \Big( \frac{\p \al}{\p t_1}(s, \, t_1, \, t_2)
  \ot \frac{\p \al}{\p t_2}(s, \, t_1, \, t_2)\Big) \, dt_1 \, dt_2 = \\
& \int_0^1 \int_0^1 \sum_{i,\, j =1}^n a_{ij}(\al (s, \, t_1, \, t_2)) \, 
  \frac{\p \al^i}{\p t_1}(s, \, t_1,\, t_2) \, 
  \frac{\p \al^j}{\p t_2}(s, \, t_1,\, t_2) \, dt_1 \, dt_2 \, , 
\end{split}
\end{equation}
where $\al^i (s, \, t_1, \, t_2) = x^i ( \al (s, \, t_1, \, t_2)) \in
\br$.  Likewise, set
$$  \gamma^i = x^i ( \gamma (t_1, \, t_2)) \in \br .  $$
\begin{lemma} \label{3.6}
A necessary condition that $J(\gamma)$ in equation \eqref{3.4} be an
extreme value for the variation $J( \al(s))$ is that
\begin{align*}
& \sum_{i,\, j =1}^n \Big( \frac{-\p a_{ij}}{\p x^{\ell}} + 
  \frac{\p a_{{\ell}j}}{\p x^i} + \frac{\p a_{i{\ell}}}{\p x^j}\Big)
  (\gamma(t_1, \, t_2)) \, \frac{\p \gamma^i}{\p t_1} \, \frac{\p
  \gamma^j}{\p t_2} \\
& + \sum_{j=1}^n (a_{\ell j} + a_{j \ell})( \gamma(t_1, \, t_2)) \, 
  \frac{\p^2 \gamma^j}{\p t_1 \, \p t_2} = 0
\end{align*}
for each $\ell = 1, \, 2, \, 3, \, \ld \, n$.  
\end{lemma} 
\begin{proof}
One computes $\frac{d( J \al (s))}{ds}$ directly and equates
$$  \frac{d(J \al (s))}{ds} \Bigr\rvert_{s=0} = 0 .  $$
In the following $a_{ij}$ and 
$\frac{ \p a_{ij}}{\p x^{\ell}}$ are evaluated at
$\al(s, \, t_1, \, t_2)$ while all partial derivatives of the $\al^i$'s
are evaluated at $(s, \, t_1, \, t_2)$.  Then
\begin{align*}
& \frac{d(J \al (s))}{ds} =   \\
& \int_0^1 \int_0^1 \Bigg\{ \sum_{i,\, j =1}^{n} \Big( \sum_{\ell =1}^{n} 
 \frac{\p a_{ij}}{\p x^{\ell}} \, \frac{\p \al^{\ell}}{\p s} \, 
 \frac{\p \al^i}{\p t_1} \, \frac{\p \al^j}{\p t_2}\Big)   \\
& \ \ \ \ \ \ \ \ \ \ \ 
  + \sum_{i,\, j =1}^{n} a_{ij} \, \frac{\p^2 \al^i}{\p s \, \p t_1}
	\, \frac{\p \al^j}{\p t_2} \, + \, \sum_{i, \, j =1}^{n} a_{ij}
	\, \frac{\p \al^i}{\p t_1} \, \frac{\p^2 \al^j}{\p s \p t_2} \Bigg\}
	\, dt_1 \, dt_2 \, ,  
\end{align*}
which can be simplified using intergration by parts and, in certain
terms, the boundary values of $\al$.  For example,
\begin{align*}
& \int_0^1 \int_0^1 a_{ij} \, \frac{\p^2 \al^i}{\p s \, \p t_1} \, 
	\frac{\p \al^j}{\p t_2} \, dt_1 \, dt_2 
  \, = \, ({\text I}) + ({\text{II}}),  \\
& ({\text I}) \ = \int_0^1 \Big[ \, a_{ij} \, \frac{\p \al^i}{\p s} \, 
	\frac{\p \al^j}{\p t_2} \, \Bigr\rvert_{t_1 = 0}^{t_1 = 1} \,\,
	dt_2 = 0, \\
& ({\text{II}})\  = - \int_0^1 \int_0^1 \Big\{ \Big( \sum_{k=1}^n
  \frac{\p a_{ij}}{\p x^k} \, \frac{\p \al^k}{\p t_1} \, \frac{\p
  \al^i}{\p s} \, \frac{\p \al^j}{\p t_2}\Big) + a_{ij} \,
  \frac{\p \al^i}{\p s}
  \, \frac{\p^2 \al^j}{\p t_1 \, \p t_2} \Big\} \, dt_1 \,
  dt_2  .   
\end{align*}   
After reindexing,
\begin{align*}
& \frac{d( J \al(s))}{ds}\Bigr\rvert_{s=0} =  \\
& \sum_{\ell =1}^n \int_0^1 \int_0^1 \frac{\p \al^{\ell}}{\p s} \, 
  \Bigg\{ \sum_{i,\, j =1}^n \Big( \frac{\p a_{ij}}{\p x^{\ell}} - 
  \frac{\p a_{\ell j}}{\p x^i} - \frac{\p a_{i \ell}}{\p x^j}\Big) \,
  \frac{\p \gamma^i}{\p t_1} \, \frac{\p \gamma^j}{\p t_2} \\
&  \hskip2in - \sum_{j=1}^n \, (a_{\ell j} + a_{j \ell}) \, 
   \frac{\p^2 \gamma^j}{\p t_1 \, \p t_2}\Bigg\} \, dt_1 \, dt_2 , 
\end{align*}
where $a_{ij}$ and $\frac{\p a_{ij}}{\p x^{\ell}}$ are evaluated at
$\gamma (t_1, \, t_2)$ and $\frac{\p \al^{\ell}}{\p s}$ is evaluated
at $(0, \, t_1, \, t_2)$.  The lemma now follows from the standard
techniques of the calculus of variations, for example \cite[p.\
432--438]{SpivakI}, and in particular \cite[p.\ 435]{SpivakI}.  
\end{proof}

Recall that the symbols 
$$ \frac{\p}{\p x^1}, \ \frac{\p}{\p x^2}, \ \ld , \ \frac{\p}{\p x^n}$$
may be interpreted as vector fields on $U$ as well as derivations of
the ring $C^{\infty}(U)$.  To state the next lemma, we introduce the
composition operators
$$  \frac{\p}{\p x^{\ell}} \circ dx^p : \chi (U) \to C^{\infty}(U)$$ 
given by
$$ \Big( \frac{\p}{\p x^{\ell}} \circ dx^p \Big)
	\Big( \sum_{i=1}^n a_i \, \frac{\p}{\p x^i}\Big) =
	\frac{\p a_p}{\p x^{\ell}}\, .  $$
\begin{lemma}\label{3.7}
Let $\om$ be a two-tensor on $M$ with local expression on $U$ 
$$  \om = \sum_{p, \, q=1}^n a_{pq} \, dx^p \ot dx^q , $$
where each $a_{pq}: U \to \br$ is $C^{\infty}$.  Then
\begin{align*}
& d \om = \\
& \sum_{p, \, q =1}^n \Bigg\{ \sum_{\ell =1}^n \frac{\p a_{pq}}{\p x^{\ell}}
  \, dx^{\ell} \ot dx^p \ot dx^q 
 - \sum_{\ell =1}^n \frac{\p a_{pq}}{\p x^{\ell}} \, dx^p \ot
   dx^{\ell} \ot dx^q  \\
& \hskip.3in + \sum_{\ell =1}^n \frac{\p a_{pq}}{\p x^{\ell}} \, dx^p \ot
   dx^q \ot dx^{\ell}  \\
& + \sum_{\ell =1}^n a_{pq} \, dx^{\ell} \ot dx^p \ot \Big( \frac{\p}{\p x^{\ell}}
  \circ dx^q \Big) + 
  \sum_{\ell =1}^n a_{pq} \, dx^{\ell} \ot dx^q \ot \Big( \frac{\p}{\p x^{\ell}}
  \circ dx^p \Big) \Bigg\} .
\end{align*}
\end{lemma}
\begin{proof}
Let $X_1, \ X_2, \ X_3 \in \chi(M)$ with local expressions
$$ X_1 = \sum_{i_1 = 1}^n c_{i_1} \, \frac{\p}{\p x^{i_1}}, \ \  
   X_2 = \sum_{i_2 = 1}^n c_{i_2} \, \frac{\p}{\p x^{i_2}}, \ \ 
   X_3 = \sum_{i_3 = 1}^n c_{i_3} \, \frac{\p}{\p x^{i_3}} \, .  $$
Then 
$$ d \om (X_1 \ot X_2 \ot X_3) =
   \sum_{i_1, \, i_2, \, i_3 =1}^n d\om \Big( c_{i_1} \frac{\p}{\p x^{i_1}}
   \ot c_{i_2} \frac{\p}{\p x^{i_2}}\ot c_{i_3} \frac{\p}{\p x^{i_3}}
   \Big) ,  $$
and from equation \eqref{2.1}
\begin{align*}
& d\om \Big( c_{i_1} \frac{\p}{\p x^{i_1}}
  \ot c_{i_2} \frac{\p}{\p x^{i_2}}\ot c_{i_3} \frac{\p}{\p x^{i_3}}\Big)= \\
& c_{i_1}c_{i_2}c_{i_3} \, \Big( \frac{\p}{\p x^{i_1}} (a_{i_2 i_3})
  - \frac{\p}{\p x^{i_2}} (a_{i_1 i_3})+\frac{\p}{\p x^{i_3}} (a_{i_1 i_2})
   \Big)  \\
& + (a_{i_2 i_3} + a_{i_3 i_2}) \, c_{i_1} c_{i_2} \, 
  \frac{\p c_{i_3}}{\p x^{i_1}} \, .
\end{align*}
Applying the right-hand side of $d \om$ in the statement of the lemma
to 
$$ c_{i_1} \frac{\p}{\p x^{i_1}}
   \ot c_{i_2} \frac{\p}{\p x^{i_2}}\ot c_{i_3} \frac{\p}{\p x^{i_3}},$$
the same result is obtained.
\end{proof}

Thus, given $\om = \sum_{p, \, q =1}^n a_{pq} \, dx^p \ot dx^q$, we
have
\begin{equation} \label{3.8}
d \om \Big( \frac{\p}{\p x^i} \ot \frac{\p}{\p x^{\ell}} \ot
  \frac{\p}{\p x^j} \Big) =
  \frac{\p a_{\ell j}}{\p x^i} - \frac{\p a_{ij}}{\p x^{\ell}} 
  + \frac{\p a_{i \ell}}{\p x^j} \, .  
\end{equation}
Although $d \om$ is $C^{\infty}(M)$-linear in the first two tensor
factors, this is not the case for the third factor:
\begin{equation} \label{3.9}
\begin{split}
& d \om \Big( \frac{\p}{\p x^i} \ot \frac{\p}{\p x^{\ell}} \ot
   c \, \frac{\p}{\p x^j} \Big) =  \\
& c \Big( \frac{\p a_{\ell j}}{\p x^i} - \frac{\p a_{ij}}{\p x^{\ell}} 
  + \frac{\p a_{i \ell}}{\p x^j} \Big) + 
  \big( a_{j \ell} + a_{\ell j} \Big) \frac{\p c}{\p x^i} \, . 
\end{split}
\end{equation}
The coefficients of $d \om$ appearing in equation \eqref{3.9} are the
same as those in lemma \ref{3.6}.  

\begin{lemma}  \label{3.10}
Let $M$ be a Riemannian manifold with metric tensor
$$ \om = \langle \ , \ \rangle = \sum_{p, \, q =1}^n g_{pq} \, dx^p \ot dx^q,
  \ \ g_{pq} = g_{qp},  $$
which is compatible with the Levi-Civita connection $\nabla$.  Then
\begin{align*}
& {\text{(i)}} \ \ d \om \Big( \frac{\p}{\p x^i} \ot \frac{\p}{\p x^{\ell}} \ot
   \frac{\p}{\p x^j} \Big) =  2[ij, \, \ell], \ \ {\text{(twice the
  Christoffel symbol),}} \\
& {\text{(ii)}} \ \ 
  d \om ( X \ot Y \ot Z) = 2 \langle Y, \, \nabla_X Z \rangle \ \
  {\text{for}} \ X, \ Y, \ Z \in \chi(M).
\end{align*}
\end{lemma} 
\begin{proof}
Part (i) follows from equation \eqref{3.8} and the definition of the
Christoffel symbols of the first kind
$$  [ij, \, \ell ] = \frac{1}{2} \Big(\frac{\p g_{j \ell}}{\p x^i} +
  \frac{\p g_{i \ell}}{\p x^j} - \frac{\p g_{ij}}{\p x^{\ell}} \Big).$$
For part (ii), recall that since $\nabla$ is compatible with the
metric tensor \cite[p.\ 54]{doCarmo}, we have
$$  X( \langle Y, \, Z \rangle ) = \langle \nabla_X Y, \, Z \rangle +
    \langle Y, \, \nabla_X Z \rangle.  $$
Since the Levi-Civita connection is symmetric (i.e. torsion-free)
\cite[p.\ 54--55]{doCarmo} \cite[p.\ 255--256]{SpivakII}, we have
$$ [X, \, Y] = \nabla_X (Y) - \nabla_Y (X).  $$
The lemma now follows from \eqref{2.1}, i.e.
\begin{align*}
d \om (X \ot Y \ot Z) & = X( \langle Y, \, Z\rangle) - Y( \langle X, \, Z\rangle) 
  + Z( \langle X, \, Y\rangle)   \\
& - \langle [X, \, Y], \, Z \rangle + \langle [X, \, Z], \, Y \rangle +
  \langle X, \, [Y, \, Z] \rangle .  
\end{align*}
\end{proof}

Thus, the obstruction to $C^{\infty}(M)$-linearity of the Leibniz
coboundary in lemma \ref{3.10} is the failure of the connection
$\nabla$ to be $C^{\infty}(M)$-linear in its second argument.  Also,
the first variation of arc length of $\gamma : I \to U$ can be
recovered as 
$$ \frac{1}{2} \, d \om \Big( \frac{d \gamma}{dt} \ot Y \ot 
\frac{d \gamma}{dt}\Big) = \Big\langle Y, \, \nabla_{\frac{d \gamma}{dt}}
\Big( \frac{d \gamma}{ dt }\Big) \Big\rangle  \, , $$
where $Y$ is a vector field which represents the variation of the
curve $\gamma$.  A necessary condition that $\gamma$ be a geodesic is
that 
$$ d \om \Big( \frac{d \gamma}{dt} \ot Y \ot 
\frac{d \gamma}{dt}\Big) = 0  $$
for all $Y$, since $\nabla_{\frac{d \gamma}{dt}}
\big( \frac{d \gamma}{ dt }\big)$ must vanish along such curves.  

The only two-tensors which are global cocycles, however, must be
two-forms.
\begin{lemma}
Let $\om$ be a two-tensor on $M$ with $d \om = 0$ in the Lebiniz
cochain complex.  Then $\om$ is a two-form.
\end{lemma}
\begin{proof}
Letting $\om = \sum_{p, q = 1}^n a_{pq} \, dx^p \ot dx^q$ be a local
expression for $\om$, then
$$  d \om \Big( \frac{\p }{\p x^i} \ot \frac{\p}{\p x^{\ell}} \ot
    \frac{\p}{\p x^j}\Big) = 0  $$
implies $\frac{\p a_{\ell j}}{\p x^i} - \frac{\p a_{ij}}{\p x^{\ell}}
+ \frac{ \p a_{i \ell}}{\p x^j} = 0$.  Furthermore
$$  d \om \Big( \frac{\p }{\p x^i} \ot \frac{\p}{\p x^{\ell}} \ot
    x^i \, \frac{\p}{\p x^j}\Big) = 0  $$
and \eqref{3.9} imply that $a_{j \ell} = - a_{\ell j}$.  Thus, $ \om $
is skew-symmetric.  
\end{proof}

The Leibniz coboundary of a $k$-tensor $ \om $ agrees with the terms
occurring in the optimization of the integral of $ \om $ if $ \om $ is
skew-symmetric in its last $(k-1)$-arguments.  Such tensors naturally
occur in the $E^2$ term of the Pirashvili spectral sequence for
Leibniz cohomology \cite{Pirashvili} \cite{Lodder}.  Suppose that $\g$
is a Lie algebra over $\br$, 
$$ \g' =  {\text{Hom}}( \g , \, \br )  $$
the coadjoint representation of $\g$, and that 
$$  H_{\text{Lie}}^* ( \g; \, \g')  $$
denotes the Lie algebra cohomology of $\g$ with coefficients in $\g'$.
An element of $H_{\text{Lie}}^{k-1}(\g; \, \g')$ 
can be represented by a tensor
$$  \al : \g^{\ot k} \to \br  $$
which is skew-symmetric in its last $(k-1)$-tensor factors
\cite{Lodder}.  Moreover, $H_{\text{Lie}}^* ( \g; \, \g')$  occurs
in the $E^2$ term of the Pirashvili spectral sequence for $HL^* (\g )$.

\section{The Local Coboundary Formula}

Let $\om$ be a $k$-tensor on $M$ with local expression
\begin{equation} \label{4.1}
\sum_{I} a_I \, dx^{i_1} \ot dx^{i_2} \ot \ \ld \ \ot dx^{i_k},
\end{equation}
where $I$ is the multi-index $(i_1, \, i_2, \, \ld , \, i_k)$, and the
summation ranges over
$$  0 \leq i_1 \leq n, \ \ 0 \leq i_2 \leq n, \ \ \ld ,\ \  0 \leq i_k \leq n.$$
To state the coboundary formula, the following local operators are
introduced.  Let
\begin{align*}
L(\om ) = & \sum_I \Bigg\{ \sum_{\ell =1}^n \frac{\p a_I}{\p x^{\ell}} \,
        dx^{\ell} \ot dx^{i_1} \ot dx^{i_2} \ot \ \ld \ \ot dx^{i_k} \\
	& - \sum_{\ell =1}^n \frac{\p a_I}{\p x^{\ell}} \,
	dx^{i_1} \ot dx^{\ell} \ot dx^{i_2} \ot \ \ld \ \ot dx^{i_k} \\
	& + \sum_{\ell =1}^n \frac{\p a_I}{\p x^{\ell}} \,
	dx^{i_1} \ot dx^{i_2} \ot dx^{\ell} \ot \ \ld \ \ot dx^{i_k} \\
	& + \ \ld  \\
	& + (-1)^{k+2} \sum_{\ell =1}^n \frac{\p a_I}{\p x^{\ell}} \,
  dx^{i_1} \ot dx^{i_2} \ot \ \ld \ \ot dx^{i_k} \ot dx^{\ell} \Bigg\},
\end{align*}
and  
\begin{align*}
& S(dx^{i_1} \ot dx^{i_2} \ot \ \ld \ \ot \ dx^{i_k}) =  \\
& \sum_{\ell =1}^n dx^{\ell} \ot dx^{i_1} \ot \Big( \frac{\p}{\p x^{\ell}}
  \circ dx^{i_2}\Big) \ot dx^{i_3} \ot \ \ld \ \ot dx^{i_k} \\
& + \sum_{\ell =1}^n dx^{\ell} \ot dx^{i_1} \ot dx^{i_2} \ot 
	\Big( \frac{\p}{\p x^{\ell}}
  \circ dx^{i_3} \Big) \ot \ \ld \ \ot dx^{i_k} \\
& + \ \ld  \\
& + \sum_{\ell =1}^n dx^{\ell} \ot dx^{i_1} \ot dx^{i_2} \ot \ \ld \ 
	\ot \ dx^{i_{k-1}}
  \ot \Big( \frac{\p}{\p x^{\ell}} \circ dx^{i_k}\Big) \\
& +(-1)^4 \sum_{\ell =1}^n dx^{\ell} \ot dx^{i_2} \ot 
	\Big( \frac{\p}{\p x^{\ell}}
  \circ dx^{i_1} \Big) \ot dx^{i_3} \ot \ \ld \ \ot dx^{i_k} \\ 
& +(-1)^5 \sum_{\ell =1}^n dx^{\ell} \ot dx^{i_2} \ot dx^{i_3} \ot
  \Big( \frac{\p}{\p x^{\ell}} \circ dx^{i_1} \Big) 
  \ot dx^{i_4} \ot \ \ld \ \ot dx^{i_k} \\ 
& + \ \ld  \\
& +(-1)^{k+2} \sum_{\ell =1}^n dx^{\ell} \ot dx^{i_2} \ot dx^{i_3} \ot
  \ \ld \  \ot dx^{i_k} \ot
  \Big( \frac{\p}{\p x^{\ell}} \circ dx^{i_1} \Big) .  
\end{align*}
\begin{theorem}
If $\om$ is a $k$-tensor on $M$, $k \geq 2$, with local expression
given in equation \eqref{4.1}, then locally the Leibniz coboundary of $\om$ is
\begin{align*}
d(\om ) = L(\om ) + & \sum_I a_I \, \Big\{ S(dx^{i_1} \ot dx^{i_2} \ot \
	\ld \ \ot dx^{i_k})  \\
 & - dx^{i_1} \ot S(dx^{i_2} \ot dx^{i_3} \ot \ \ld \ \ot dx^{i_k})\\
 & + dx^{i_1} \ot dx^{i_2} \ot S(dx^{i_3} \ot \ \ld \ \ot dx^{i_k})\\
 & + \ \ld  \\
 & + (-1)^{k-2} dx^{i_1} \ot \ \ld \ \ot dx^{i_{k-2}} \ot 
	S(dx^{i_{k-1}} \ot dx^{i_k}) \Big\} .  
\end{align*}
\end{theorem}
\begin{proof}
The proof proceeds by induction on $k$ with the case $k=2$ proven in
lemma \ref{3.7}.  To streamline the inductive step, define a
$(k-1)$-tensor $\beta$ (often called a contraction) by
$$  \beta(v_2 \ot \, \ld \, \ot v_k ) = 
   \om \Big( \frac{\p}{\p x^{j_1}} \ot v_2 \ot \, \ld \, \ot v_k \Big),$$
where $\frac{\p}{\p x^{j_1}}$ is a fixed canonical vector field on a
coordinate chart.  Let
$$ z = \Big( c_{j_1}\frac{\p}{\p x^{j_1}} \ot c_{j_2}\frac{\p}{\p x^{j_2}} \ot
   c_{j_3}\frac{\p}{\p x^{j_3}} \ot \ \ld \ \ot
   c_{j_{k+1}}\frac{\p}{\p x^{j_{k+1}}} \Big)  $$
Then
\begin{align*}
& d \om (z) =  \\
& = c_{j_1}c_{j_2} \, \ld \, c_{j_{k+1}} \, \frac{\p}{\p x^{j_1}}
    \big( a_{j_2 \, j_3 \, \ld \, j_{k+1}} \big) \\
& + \sum_{m=3}^{k+1} \big( a_{j_2 \, j_3 \, \ld \, j_{k+1}} + 
 (-1)^{m+1} a_{j_m \, j_2 \, j_3 \, \ld \, {\hat{j}}_m \, \ld \, j_{k+1}}\big)
 \, c_{j_1} c_{j_2} \, \ld \, \hat{c}_{j_m} \, \ld \, c_{j_{k+1}} \, 
 \frac{\p c_{j_m}}{\p x^{j_1}} \\
& + \sum_{i_1 = 1}^n (- dx^{i_1} \ot d \beta) (z)  \\
& = \sum_I \sum_{\ell =1}^n \Big( \frac{\p a_I}{\p x^{\ell}} \, dx^{\ell}
 \ot dx^{i_1} \ot dx^{i_2} \ot \ \ld \ \ot dx^{i_k} \Big) (z)  \\
& + \sum_I a_I \, S ( dx^{i_1} \ot dx^{i_2} \ot \ \ld \ \ot dx^{i_k})(z)
  + \sum_{i_1 =1}^n - (dx^{i_1} \ot d \beta) (z) ,
\end{align*}
whence follows the result.  
\end{proof}

\begin{lemma}
If $\om$ is a $k$-tensor on an $n$-dimensional differentiable manifold
$M$, $k \leq (n+1)$, and $d \om = 0$, then $\om$ is a $k$-form.
\end{lemma}
\begin{proof}
Clearly $d \om (\Xi ) = 0$ for any $\Xi \in \chi(M)^{\ot (k+1)}$.  Since
$$  d \om \Big( \frac{\p}{\p x^{j_1}} \ot \frac{\p}{\p x^{j_2}} \ot
    \ \ld \ \ot \frac{\p}{\p x^{j_{k+1}}}\Big) = 0 ,  $$
we have (using the notation in equation \eqref{4.1}), 
\begin{align*}
& \frac{\p}{\p x^{j_1}} \big( a_{j_2 \, j_3 \, \ld \, j_{k+1}}\big) 
 -\frac{\p}{\p x^{j_2}} \big( a_{j_1 \, j_3 \, j_4 \, \ld \, j_{k+1}}\big) 
 + \ \ld  \\
& \hskip1in + \, \cdots \, + (-1)^{k+2} \frac{\p}{\p x^{j_{k+1}}}
  \big( a_{j_1 \, j_2 \, j_3 \, \ld \, j_k} \big) = 0 .  
\end{align*}
Choosing $j_1$, $j_2$, $\ld$, $j_n$ to be distinct (which may be done
since dim$(M) = n$), we also have
\begin{align*}
& d \om \Big( \frac{\p}{\p x^{j_1}} \ot \frac{\p}{\p x^{j_2}} \ot 
 \ \ld \ \ot \frac{\p}{\p x^{j_{q-1}}} \ot x^{j_p}
 \frac{\p}{\p x^{j_q}}  \\
& \hskip2in \ot \frac{\p}{\p x^{j_{q+1}}} \ot \ \ld \ \ot 
 \frac{\p}{\p x^{j_{k+1}}} \Big) = 0 , 
\end{align*}
where $p \leq q - 2$.  Thus, 
$$ a_{j_1 \, j_2 \, \ld \, \hat{j}_p \, \ld \, j_{k+1}} + (-1)^{q-p}
   a_{j_1 \, j_2 \, \ld \, j_{p-1} \, j_q \, j_{p+1} \, \ld \, 
   \hat{j}_q \, \ld \, j_{k+1}} = 0 .  $$
\end{proof}

The section is closed with the computation of the Leibniz coboundary
of one of most important tensors in differential geometry, the Riemann
curvature tensor, $R$.  Let $\nab$ be the Levi-Civita connection on a
Riemannian manifold $M$ with metric $\langle \ , \ \rangle$.  Given
$X$, $Y$, $Z$, $W$ $\in \chi(M)$, then $R$ is the four-tensor defined
by \cite{SpivakII}
$$ R(X \ot Y \ot Z \ot W) = \langle \nab_X \nab_Y (Z) - 
  \nab_Y \nab_X (Z) - \nab_{[X, \, Y]}(Z), \ W \rangle .  $$
The Leibniz coboundary $dR$ is expressed in terms of the covariant
derivative $\nab R$, which we briefly review.  Let $\om$ be a
$k$-tensor on $M$ and
$$  X_1, \ X_2, \ \ld , \ X_k, \ Z \in \chi (M) .  $$
Then \cite[p.\ 102]{doCarmo} the covariant derivative $\nab \om$ is
the $(k+1)$-tensor given by
\begin{align*}
 \nab \om & (X_1 \ot X_2 \ot \ \ld \ \ot X_k \ot Z ) = 
  Z \big( \om (X_1 \ot X_2 \ot \ \ld \ \ot X_k) \big) \\
& - \om \big( \nab_Z (X_1) \ot X_2 \ot \ \ld \ \ot X_k ) 
  - \om \big( X_1 \ot \nab_Z (X_2) \ot \ \ld \ \ot X_k \big) \\
& - \, \cdots \, - \om \big( X_1 \ot X_2 \ot \ \ld \ \ot \nab_Z (X_k)\big).
\end{align*}
Note that $\nab \om (X_1 \ot X_2 \ot \ \ld \ \ot X_k \ot Z)$ is often
denoted as $(\nab_Z \om )(X_1 \ot \ \ld \ot X_k )$.  The following
properties of $R$ are useful in the computation of $dR$, where $X$,
$Y$, $Z$, $W$, $T$ $\in \chi (M)$:

\noindent
(i) Bianchi's secnod identity \cite[p.\ 106]{doCarmo}
\begin{equation} \label{4.4}
\begin{split}
& (\nab_T R )(X \ot Y \ot Z \ot W) + (\nab_Z R )(X \ot Y \ot W \ot T) \\
& \hskip2in + 
       (\nab_W R)(X \ot Y \ot T \ot Z) = 0 ,
\end{split}
\end{equation}
which may also be expressed as \cite[p.\ 34]{Sakai}
\begin{equation} \label{4.5}
\begin{split}
& (\nab_X R )(Y \ot Z \ot W \ot T) + (\nab_Y R)(Z \ot X \ot W \ot T) \\
& \hskip2in + 
	(\nab_Z R)(X \ot Y \ot W \ot T ) = 0,  
\end{split}
\end{equation}
(ii) skew-symmetry in certain coordinates \cite[p.\ 91]{doCarmo}
\begin{equation} \label{4.6}
\begin{split}
  &  R (X \ot Y \ot Z \ot T) = - R (Y \ot X \ot Z \ot T)  \\
  &  R (X \ot Y \ot Z \ot T) = - R (X \ot Y \ot T \ot Z)  \\
  &  R (X \ot Y \ot Z \ot T) = + R (Z \ot T \ot X \ot Y) .
\end{split}
\end{equation}

\begin{lemma}\label{4.7}
Let $M$ be a Riemannian manifold with curvature tensor $R$ and
Levi-Civita connection $\nab$.  Then for 
$$ X, \ Y, \ Z, \ W, \ T \in \chi(M) ,  $$
one has
\begin{align*}
dR & (X \ot Y \ot Z \ot W \ot T) = -(\nab_Z R) (X \ot Y \ot W \ot T)\\
& + R(Z \ot T \ot Y \ot \nab_X W) - R(Y \ot Z \ot T \ot \nab_X W)  \\
& - R(Z \ot W \ot Y \ot \nab_X T) + R(Y \ot Z \ot W \ot \nab_X T)  \\
& - R(Z \ot T \ot X \ot \nab_Y W) + R(X \ot Z \ot T \ot \nab_Y W)  \\
& + R(Z \ot W \ot X \ot \nab_Y T) - R(X \ot Z \ot W \ot \nab_Y T) .
\end{align*}
\end{lemma}
\begin{proof}
After applying equations \eqref{2.1}, \eqref{4.4}, \eqref{4.6}, and
symmetry of the connection,
$$  [X, \, Y] = \nab_X Y - \nab_Y X ,  $$
one has
\begin{align*}
dR & (X \ot Y \ot Z \ot W \ot T ) =  \\
& X \big( R(W \ot T \ot Y \ot Z)\big) - Y \big( R(W \ot T \ot X \ot Z)\big)\\
& - R(W \ot T \ot \nab_X Y \ot Z) + R(W \ot T \ot \nab_Y X \ot Z)  \\
& + R(W \ot T \ot \nab_X Z \ot Y) - R(Z \ot T \ot \nab_X W \ot Y)  \\
& + R(Z \ot W \ot \nab_X T \ot Y) + R(W \ot T \ot X \ot \nab_Y Z)  \\
& - R(Z \ot T \ot X \ot \nab_Y W) + R(Z \ot W \ot X \ot \nab_Y T) .
\end{align*}
The lemma now follows form equations \eqref{4.5} and \eqref{4.6}.
\end{proof}

In the statement of lemma \ref{4.7}, note that $\nab_Z R$ is the
$C^{\infty}(M)$-linear term of $dR$, while the remaining eight terms comprise
the non-linear pieces.

\end{document}